\documentclass{article}
\usepackage{mathtools}
\usepackage{amssymb}
\usepackage{tikz}
\usepackage{caption}
\usetikzlibrary{decorations.markings, decorations.pathreplacing}
\usepackage{float}
\usepackage[margin=2.25cm]{geometry}
\usepackage{mathdots}
\usepackage{enumitem}
\usepackage{titling}
\setlist[itemize]{leftmargin=*}

\begin{document}

\title{\huge An averaging result for union-closed families of sets}
\author{\Large Christopher Bouchard}
\date{}
\maketitle

\abstract{\vspace{0.1cm} \noindent Let $\mathcal{A}$ be a union-closed family of sets with base set $b(\mathcal{A})=\bigcup_{A \in \mathcal{A}}A$ denoted by $[n]=\{1, \cdots, n\}$, and for any real $x>0$, let $\mathcal{A}_{<x} = \{A \in \mathcal{A} \ | \ |A| < x\}$. Also, denote by $\mathcal{B}$ any smallest irredundant subfamily of $\mathcal{A}_{<n/2}$ such that $b(\mathcal{B})=b(\mathcal{A}_{<n/2})$. We prove that if $\mathcal{A}$ is separating with height $h = 4 \leq n$ and $0 \leq |\mathcal{B}| \leq 2$, then the average size of a member set from $\mathcal{A}$ is at least $n/2$. We show that $h=4$ is greatest possible with respect to this result, and conclude by considering the remaining domain $3 \leq |\mathcal{B}| \leq 4$.}

\bigskip

\section*{1. Introduction}

A family of sets $\mathcal{A}$ is \textit{union-closed} if it is finite with at least one nonempty member set, all of its member sets are finite and distinct, and $X,Y \in \mathcal{A}$ implies that $X \cup Y \in \mathcal{A}$. Denote the \textit{base set} $\bigcup_{A \in \mathcal{A}} A$ of $\mathcal{A}$ by $[n]=\{1, \cdots, n\}$ (and generally of a family $\mathcal{F}$ by $b(\mathcal{F}) \coloneqq \bigcup_{F \in \mathcal{F}}F$). The union-closed sets conjecture, also referred to as Frankl's conjecture, is a well-known problem that has been studied from various perspectives, including those of graph and lattice theory (see [2] and [10], respectively, for details of the pertinent formulations), and also information theory (see [9] for the initial breakthrough work and [4] for a description of several follow-up results). The conjecture is stated as follows:

\medskip

\noindent \textbf{Conjecture 1.1.} \textit{For any union-closed family} $\mathcal{A}$\textit{, there exists an element of} $[n]$ \textit{that is in at least} $\frac{|\mathcal{A}|}{2}$ \textit{member sets of} $\mathcal{A}$\textit{.}

\medskip

For two sets $X_1$ and $X_2$, we use the notation $X_1 \subsetneq X_2$ to denote that $X_1$ is a proper subset of $X_2$, i.e. $(X_1 \subseteq X_2) \land (X_1 \neq X_2)$. A \textit{chain} $\mathcal{C}$ in a finite family of sets $\mathcal{A}$ is a subfamily of $\mathcal{A}$ such that $(X_1,X_2 \in \mathcal{C}) \ \land \ (X_1 \neq X_2) \implies (X_1 \subsetneq X_2) \ \lor \ (X_2 \subsetneq X_1)$, and the \textit{height} of $\mathcal{A}$, denoted by $h$, is the maximum size of a chain in $\mathcal{A}$. Conjecture 1.1 has been shown to hold for any union-closed family $\mathcal{A}$ with $h \leq 3$ (see [13], as well as [5]). In general, we have the following:

\medskip

\noindent \textbf{Theorem 1.2.} \textit{For any union-closed family} $\mathcal{A}$ \textit{with} $|\mathcal{A}|>1$\textit{, there exists an element of} $[n]$ \textit{that is in at least} $\frac{|\mathcal{A}|+h-3}{h-1}$ \textit{member sets of} $\mathcal{A}$\textit{.}

\medskip

\noindent \textit{Proof.} Let $\mathcal{A}$ be a union-closed family with $|\mathcal{A}|>1$, and consider any chain $\mathcal{C} = \{C_1, \cdots, C_h\}$ in $\mathcal{A}$ of maximum size, where $1 \leq i < j \leq h$ implies that $C_j \subsetneq C_i$ without loss of generality. Since $[n] \in \mathcal{A}$ and every member set of $\mathcal{A}$ is a subset of $[n]$, we have that $C_1=[n]$. For every $i \in [h-1]$, let $c_i$ be some element from $C_i \setminus C_{i+1}$. We consider any $X \in \mathcal{A} \setminus \{C_1, C_h\}$, and assume that $c_i \not \in X$ for all $i \in [h-1]$. It must be that $X \cap (C_1 \setminus C_2) = \emptyset$. (Otherwise, $\emptyset \subsetneq X \cap (C_1 \setminus C_2) \subsetneq C_1 \setminus C_2$, which implies that $\{C_1, C_2 \cup X, C_2, \cdots, C_h\}$ is a chain of size $h+1$ in $\mathcal{A}$, contradicting the definition of $h$.) Now consider any $i \in [h-1]$ such that $i \geq 2$. If, for all $j \in [i-1]$, $X \cap (C_j \setminus C_{j+1}) = \emptyset$, then $X \cap (C_{i} \setminus C_{i+1}) = \emptyset$. (Otherwise, $\emptyset \subsetneq X \cap (C_i \setminus C_{i+1}) \subsetneq C_i \setminus C_{i+1}$, which implies that $\{C_1, \cdots, C_i, C_{i+1} \cup X, C_{i+1}, \cdots, C_h\}$ is a chain of size $h+1$ in $\mathcal{A}$, again contradicting the definition of $h$.) Therefore, we have by induction that $X \cap (C_i \setminus C_{i+1}) = \emptyset$ for every $i \in [h-1]$, which implies that $X \subsetneq C_h$. We then have that $\{C_1, \cdots, C_h, X\}$ is a chain in $\mathcal{A}$ of size $h+1$, once again a contradiction with the definition of $h$. Hence, there must exist $i \in [h-1]$ such that $c_i \in X$. It follows that, for some $j \in [h-1]$, $c_j$ is in at least $\frac{|\mathcal{A}\setminus \{C_1, C_h\}|}{h-1}=\frac{|\mathcal{A}|-2}{h-1}$ member sets of $\mathcal{A} \setminus \{C_1,C_h\}$. Then because $c_j$ is in $C_1$, it must belong to at least $\frac{|\mathcal{A}|-2}{h-1}+1=\frac{|\mathcal{A}|+h-3}{h-1}$ member sets of $\mathcal{A}$, which completes the proof of Theorem 1.2.

\medskip

A chain $\mathcal{C}$ in $\mathcal{A}$ is \textit{maximal} if $\mathcal{C} \subseteq \mathcal{C}'$ implies that $\mathcal{C}=\mathcal{C}'$ for any chain $\mathcal{C}'$ in $\mathcal{A}$. We observe that Theorem 1.2 continues to hold if we replace $h$ with $r$, where $r$ is the minimum size of a maximal chain in $\mathcal{A}$, thereby verifying Conjecture 1.1 for $r \leq 3$. Now, for $h \geq 4$, Conjecture 1.1 is stronger than Theorem 1.2. (Moreover, Theorem 1.2 is not optimal for $h \geq 4$ because in [12] it is shown that there is an element from $[n]$ in at least $\frac{3-\sqrt{5}}{2}|\mathcal{A}|$ member sets of $\mathcal{A}$.) Hence, for values of $h$ greater than or equal to $4$, Conjecture 1.1 requires other proof techniques. Averaging is a common approach to Conjecture 1.1 (see Section 6 of [3] for detailed discussion, [1], [6], and [7] for example applications, and [11] for the lowest average set size of a union-closed family in terms of family size). In the present work, we consider its application to the case $h=4$. 

\pagebreak

The following definitions pertain to any finite family of sets $\mathcal{A}$ with base set $[n]$:

\smallskip

\begin{itemize}

\item $\mathcal{A}$ is \textit{separating} if, for any two distinct elements $x$ and $y$ in $[n]$, there exists $A \in \mathcal{A}$ such that $(x \in A \ \land \ y \not \in A) \ \lor \ (x \not \in A \ \land \ y \in A)$. 

\smallskip

\item For any $\mathcal{S} \subseteq \mathcal{A}$ and $S \in \mathcal{S}$, $\texttt{irr}_{\mathcal{S}}(S) \coloneqq \{s \in S \ | \ s \not \in b(\mathcal{S} \setminus \{S\})\}$. Accordingly, $\mathcal{S}$ is \textit{irreduntant} if $S \in \mathcal{S}$ implies that $\texttt{irr}_{\mathcal{S}}(S) \neq \emptyset$.

\smallskip

\item For any real $x \geq 0$, $\mathcal{A}_{< x} \coloneqq \{A \in \mathcal{A} \ | \ |A| < x\}$, $\mathcal{A}_{\leq x} \coloneqq \{A \in \mathcal{A} \ | \ |A| \leq x\}$, $\mathcal{A}_{> x} \coloneqq \{A \in \mathcal{A} \ | \ |A| > x\}$, and $\mathcal{A}_{\geq x} \coloneqq \{A \in \mathcal{A} \ | \ |A| \geq x\}$. Similarly, for any $X \subseteq [n]$, $\mathcal{A}_{\subsetneq X} \coloneqq \{A \in \mathcal{A} \ | \ A \subsetneq X\}$ and $\mathcal{A}_{\subseteq X} \coloneqq \{A \in \mathcal{A} \ | \ A \subseteq X\}$.

\smallskip

\item We set $B = b(\mathcal{A}_{<n/2})$, and denote by $\mathcal{B} \coloneqq \mathcal{B}(\mathcal{A})$ any irredundant subfamily of $\mathcal{A}_{<n/2}$ of minimum size such that $b(\mathcal{B})=B$.

\end{itemize}

\smallskip
\smallskip

The main result of this work, proved in the next section, is that the average size of a member set in any separating union-closed family $\mathcal{A}$ with $h = 4 \leq n$ and $0 \leq |\mathcal{B}| \leq 2$ is at least $\frac{n}{2}$, which implies that Conjecture 1.1 holds for these particular families. We also demonstrate via an appropriate construction that, with respect to this averaging result, $h=4$ is greatest possible, and we conclude with partial results for when $3 \leq |\mathcal{B}| \leq 4$. The following lemma is essential to this study.

\medskip

\noindent \textbf{Lemma 1.3.} \textit{For any separating union-closed family} $\mathcal{A}$\textit{, every maximal chain} $\mathcal{C}$ \textit{in} $\mathcal{A}$ \textit{contains a member set of size} $n-1$\textit{.}

\medskip

\noindent \textit{Proof:} Let $\mathcal{A}$ be a separating union-closed family, and $\mathcal{C} = \{C_1, \cdots, C_{|\mathcal{C}|}\}$ be any maximal chain in $\mathcal{A}$, where $1 \leq i < j \leq |\mathcal{C}|$ implies that $C_j \subsetneq C_i$ without loss of generality. We have that $C_1=[n]$, as $[n] \in \mathcal{A}$ and every member set of $\mathcal{A}$ is a subset of $[n]$. Assume that $\mathcal{C}$ has no member set of size $n-1$, which is equivalent to the assumption that $|C_2|<n-1$. This implies that $|C_1 \setminus C_2| \geq 2$. Then, because $\mathcal{A}$ is separating, there must exist $A \in \mathcal{A}$ and two elements $x,y \in C_1 \setminus C_2$ such that $x \in A$ and $y \not \in A$. It follows that $\mathcal{C}'=\{C_1, C_2 \cup A, C_2, \cdots, C_{|\mathcal{C}|}\}$ is a chain in $\mathcal{A}$ such that $\mathcal{C} \subsetneq \mathcal{C}'$, which contradicts the maximality of $\mathcal{C}$. Thus, $|C_2|=n-1$, proving Lemma 1.3.

\medskip

Lemma 1.3 of course applies to any chain of maximum size $h$ in $\mathcal{A}$. Throughout this writing, we denote by $Y$ some member set of $\mathcal{A}$ such that $|Y| = n-1$, with existence guaranteed by the lemma. We also denote by $\textrm{Avg}(\mathcal{A})$ the average size of a member set in $\mathcal{A}$, i.e. $\textrm{Avg}(\mathcal{A})=\frac{\sum_{A \in \mathcal{A}} |A|}{|\mathcal{A}|}$. Considering Lemma 1.3, we conclude the current section with the following theorem:

\medskip

\noindent \textbf{Theorem 1.4.} \textit{If} $\mathcal{A}$ \textit{is a separating union-closed family with} $h \leq 3$\textit{, then} $\textrm{Avg}(\mathcal{A}) \geq \frac{n}{2}$\textit{.}

\medskip

\noindent \textit{Proof.} First, we assume that $h \leq 2$. By Lemma 1.3, $\mathcal{A}=[n] \cup \mathcal{M}$, where $M \in \mathcal{M}$ implies that $|M|=n-1$. Thus, $\textrm{Avg}(\mathcal{A}) = \frac{n + k(n-1)}{k+1}$ for some $k \in \{0,1, \cdots, n\}$, which implies that $\textrm{Avg}(\mathcal{A})\geq \frac{n}{2}$. Next, we assume that $h=3$. If $|\mathcal{A}_{<n/2}|=0$, then $|\mathcal{A}| \geq \frac{n}{2}$ for every $A \in \mathcal{A}$, implying that $\textrm{Avg}(\mathcal{A}) \geq \frac{n}{2}$. If $|\mathcal{A}_{<n/2}| = 1$, then the subfamily $\hat{\mathcal{A}}=\mathcal{A}_{<n/2}\cup\{[n]\}$ of $\mathcal{A}$ has $\textrm{Avg}(\hat{\mathcal{A}}) \geq \frac{n}{2}$, implying again that $\textrm{Avg}(\mathcal{A}) \geq \frac{n}{2}$, as any member set in $\mathcal{A} \setminus \hat{\mathcal{A}}$ has size greater than or equal to $\frac{n}{2}$. Now, if $|\mathcal{A}_{<n/2}| \geq 2$, then any distinct member sets $X_1$ and $X_2$ in $\mathcal{A}_{<n/2}$ have $X_1 \cap X_2 = \emptyset$ and $|X_1|=|X_2|=\frac{n-1}{2}$. (Otherwise, $|X_1 \cup X_2| < n-1$ and $(X_1 \subsetneq X_1 \cup X_2 \subsetneq [n]) \lor (X_2 \subsetneq X_1 \cup X_2 \subsetneq [n])$, contradicting Lemma 1.3.) If there exist only two member sets $X_1$ and $X_2$ in $\mathcal{A}_{<n/2}$, then the subfamily $\hat{\mathcal{A}}=\{X_1,X_2, [n]\}$ of $\mathcal{A}$ has $\textrm{Avg}(\hat{\mathcal{A}})=\frac{2n-1}{3} \geq \frac{n}{2}$, again implying that $\textrm{Avg}(\mathcal{A}) \geq \frac{n}{2}$. Else, $\mathcal{A}=\{\{1,2,3\},\{1,2\},\{1,3\},\{2,3\},\{1\},\{2\},\{3\}\}$ and $\textrm{Avg}(\mathcal{A})=\frac{12}{7} \geq \frac{n}{2}=\frac{3}{2}$, completing the proof of Theorem 1.4.

\section*{2. The main result}

We first note that $0 \leq |\mathcal{B}| \leq h$ for any union-closed family $\mathcal{A}$. (Otherwise, for such a family $\mathcal{A}$, there exists $\{B_1, \cdots, B_{h+1}\} \subseteq \mathcal{B}$, and $B_1 \subsetneq B_1 \cup B_2 \subsetneq \cdots \subsetneq \bigcup_{j \in [i]} B_j \subsetneq \cdots \subsetneq \bigcup_{j \in [h]}B_j \subsetneq \bigcup_{j \in [h+1]}B_j$, forming a chain of size $h+1$ in $\mathcal{A}$, which contradicts the definition of $h$.) Thus, for any union-closed family $\mathcal{A}$ with $h=4$, $|\mathcal{B}|$ must belong to the set $\{0,1,2,3,4\}$. We now state the main theorem of this work:

\vspace{0.5cm}

\noindent \textbf{Theorem 2.1.} \textit{For any separating union-closed family} $\mathcal{A}$ \textit{with} $h = 4 \leq n$ \textit{and} $0 \leq |\mathcal{B}| \leq 2$:

\medskip

\[\textrm{Avg}(\mathcal{A})=\frac{\sum_{A \in \mathcal{A}}|A|}{|\mathcal{A}|} \geq \frac{n}{2}\textrm{.}\]

\pagebreak

\noindent \textbf{\normalsize{Proof of Theorem 2.1}}

\medskip

\noindent The need for Theorem 2.1 to state that $4 \leq n$ is a consequence of the separating union-closed family $\mathcal{A}=\{\{1,2,3\}, \{1,2\}, \{1\}, \{2\}, \emptyset\}$ (or any relabeling thereof), which has $h=4$ and $0 \leq |\mathcal{B}| \leq 2$, yet also has $n=3$ and $\textrm{Avg}(\mathcal{A})=\frac{7}{5}<\frac{n}{2}=\frac{3}{2}$. 

\medskip

\noindent Let $\mathcal{A}$ be any separating union-closed family with $h = 4 \leq n$ and $0 \leq |\mathcal{B}| \leq 2$. For the proof of Theorem 2.1, we must show that $\textrm{Avg}(\mathcal{A}) \geq \frac{n}{2}$.

\medskip

\noindent We first establish Propositions A, B, and C, valid for $|B| < n-1$.

\medskip
\smallskip

\noindent \textbf{Proposition A.} If $X_1$ and $X_2$ are distinct member sets of $\mathcal{A}_{\subsetneq B}$, then $(B \setminus X_1) \cap (B \setminus X_2) = \emptyset$.

\medskip

\noindent \textit{Proof.} If not, then there exist some distinct $X_1,X_2 \in \mathcal{A}_{\subsetneq B}$ such that $(B \setminus X_1) \cap (B \setminus X_2) \neq \emptyset$, and we have that $(X_1 \subsetneq X_1 \cup X_2 \subsetneq B \subsetneq [n]) \lor (X_2 \subsetneq X_1 \cup X_2 \subsetneq B \subsetneq [n])$, contradicting Lemma 1.3.

\medskip
\smallskip

\noindent \textbf{Proposition B.} $(\textrm{Avg}(\mathcal{A}) \geq \frac{n}{2}) \ \lor \ (1 \leq |\mathcal{A}_{\subsetneq B}| \leq |B|)$.

\medskip

\noindent \textit{Proof.} If $|\mathcal{A}_{\subsetneq B}|<1$, then $\mathcal{A}_{\subsetneq B} = \emptyset$ and the subfamily $\hat{\mathcal{A}}=\mathcal{A}_{\subseteq B} \cup \{[n]\}$ of $\mathcal{A}$ has $\textrm{Avg}(\hat{\mathcal{A}}) \geq \frac{n}{2}$. Then, because any member set of $\mathcal{A} \setminus \hat{\mathcal{A}}$ has size greater than or equal to $\frac{n}{2}$, we have that $\textrm{Avg}(\mathcal{A}) \geq \frac{n}{2}$. Next, if $|\mathcal{A}_{\subsetneq B}|>|B|$, then $\mathcal{A}_{\subsetneq B}=\{X_1, \cdots, X_k\}$ for some $k>|B|$. By Proposition A, $(B \setminus X_1), \cdots, (B \setminus X_k)$ are all mutually disjoint. It follows that $|B| < \ |\bigcup_{i \in [k]}(B \setminus X_i)|$, which contradicts that $\bigcup_{i \in [k]}(B \setminus X_i) \subseteq B$.

\medskip
\smallskip

\noindent \textbf{Proposition C.} $\sum_{X \in \mathcal{A}_{\subsetneq B}}|X| \geq (|\mathcal{A}_{\subsetneq B}|-1)|B|$.

\medskip

\noindent \textit{Proof.} By Proposition A, we have that $\sum_{X \in \mathcal{A}_{\subsetneq B}}|B \setminus X| \leq |B|$. Proposition C then follows from the identity $\sum_{X \in \mathcal{A}_{\subsetneq B}}|X|=|B||\mathcal{A}_{\subsetneq B}|-\sum_{X \in \mathcal{A}_{\subsetneq B}}|B \setminus X|$.

\bigskip
\smallskip

\noindent We now divide the proof of Theorem 2.1 into cases based on the value of $|\mathcal{B}|$.

\bigskip
\smallskip

\noindent \textbf{\normalsize{Case 0: $(|\mathcal{B}| = 0)$}}

\medskip

\noindent Here, if $\mathcal{A}_{<n/2}=\emptyset$, then $|A| \geq \frac{n}{2}$ for every $A \in \mathcal{A}$, making $\textrm{Avg}(\mathcal{A})\geq \frac{n}{2}$. Else, if $\mathcal{A}_{<n/2}=\{\emptyset\}$, then consider the subfamily $\hat{\mathcal{A}}=\{\emptyset, [n]\}$ of $\mathcal{A}$. Because any member set in $\mathcal{A} \setminus \hat{\mathcal{A}}$ must have size greater than or equal to $\frac{n}{2}$, it follows from $\textrm{Avg}(\hat{\mathcal{A}})=\frac{n}{2}$ that $\textrm{Avg}(\mathcal{A}) \geq \frac{n}{2}$.

\bigskip
\smallskip

\noindent \textbf{\normalsize{Case 1: $(|\mathcal{B}| = 1)$}}

\medskip

\noindent In this case, $B$ is the unique member set of $\mathcal{B}$, and we have that $\mathcal{A}_{<|B|}=\mathcal{A}_{\subsetneq B}$ and $\mathcal{A}_{\leq|B|}=\mathcal{A}_{\subseteq B}$.

\medskip

\noindent We consider any $A \in \mathcal{A}_{>|B|}$. If $|A| < \frac{n}{2}$, then either $B \subsetneq A$, which contradicts that $|B|$ is greatest among member sets of $\mathcal{A}_{<n/2}$, or $B \setminus A \neq \emptyset$, implying that $|\mathcal{B}|>1$, again a contradiction.

\medskip

\noindent By Proposition B, we have that $|\mathcal{A}_{\subseteq B}| \leq |B|+1$. (We may assume that $n>4$, as $n=4$ implies that $|\mathcal{A}_{<n/2}| = 2$, which when coupled with $\{Y,[n]\} \subseteq \mathcal{A}_{>n/2}$ in turn implies that $\textrm{Avg}(\mathcal{A}) \geq \frac{n}{2}$.) It follows that $|\mathcal{A}_{> |B|} \setminus \{Y, [n]\}| \geq |\mathcal{A}|-|B|-3 \geq |\mathcal{A}|-\frac{n+5}{2} \geq |\mathcal{A}|-n$. We let $\hat{\mathcal{A}}=\mathcal{A} \setminus \tilde{\mathcal{A}}$, where $\tilde{\mathcal{A}} \subseteq \mathcal{A}_{>|B|} \setminus \{Y, [n]\}$ such that $|\tilde{\mathcal{A}}|=|\mathcal{A}|-n$, so we have that $|\hat{\mathcal{A}}| = n$. Existence of such a family $\tilde{\mathcal{A}}$ is guaranteed by the following lemma:

\medskip

\noindent \textbf{Lemma 2.1.1} (Falgas-Ravry [8])\textbf{.} \textit{If} $\mathcal{A}$ \textit{is a separating union-closed family, then} $|\mathcal{A}| \geq n$\textit{.}

\medskip

\noindent \textit{Proof:} This follows from Lemma 2 of [8]. We provide a proof by induction on family size. For any family of sets $\mathcal{F}$ and $x \in b(\mathcal{F})$, let $\mathcal{F}_{\{x\}}=\{F \in \mathcal{F} \ | \ x \in F\}$. Now, let $\mathcal{A}$ be any separating union-closed family of sets. If $|\mathcal{A}|=1$, then $|b(\mathcal{A})|=1$ and $|\mathcal{A}| \geq |b(\mathcal{A})|$. If $|\mathcal{A}| > 1$, then consider any $x \in [n]$ such that $|\mathcal{A}_{\{x\}}|=\max_{y \in [n]}\{|\mathcal{A}_{\{y\}}|\}$, and let $\hat{\mathcal{A}}_{\{x\}}=\{A \setminus \{x\} \ | \ A \in \mathcal{A}_{\{x\}}\}$. If $\mathcal{A}_{\{x\}} \subsetneq \mathcal{A}$, then let $\mathcal{A}'=\mathcal{A}_{\{x\}}$. Else, if $\mathcal{A}_{\{x\}} = \mathcal{A}$, then let $\mathcal{A}'=(\hat{\mathcal{A}}_{\{x\}})_{\{y\}}$ for some $y \in b(\hat{\mathcal{A}}_{\{x\}})$ such that $|(\hat{\mathcal{A}}_{\{x\}})_{\{y\}}|=\max_{z \in b(\hat{\mathcal{A}}_{\{x\}})}\{|{(\hat{\mathcal{A}}_{\{x\}}})_{\{z\}}|\}$. In either case, $\mathcal{A}'$ is union-closed and separating with $|\mathcal{A}'|<|\mathcal{A}|$. We have that $|\mathcal{A}'| \geq |b(\mathcal{A}')|$ by the induction hypothesis. Further, $n-1 \leq |b(\mathcal{A}')|$. It follows that $|\mathcal{A}| \geq n$, completing the proof of Lemma 2.1.1.

\medskip

\noindent Since any member set of $\tilde{\mathcal{A}}$ has size greater than or equal to $\frac{n}{2}$, it is sufficient for resolving Case 1 to prove that $\textrm{Avg}(\hat{\mathcal{A}}) \geq \frac{n}{2}$.

\pagebreak

\noindent In this regard, we introduce the lower bound $\zeta(|B|, |\mathcal{A}_{\subsetneq B}|)$ of $\textrm{Avg}(\hat{\mathcal{A}})$:

\medskip

\[\hspace{-0.675cm} \textrm{Avg}(\hat{\mathcal{A}})=\frac{\sum_{A \in \hat{\mathcal{A}}}|A|}{|\hat{\mathcal{A}}|} = \frac{|Y|+n+\sum_{A \in \mathcal{A}_{\leq |B|}}|A|+\sum_{A \in \hat{\mathcal{A}}_{>|B|} \setminus \{Y,[n]\}}|A|}{n}\]

\smallskip

\[\hspace{1.325cm} \geq \zeta(|B|, |\mathcal{A}_{\subsetneq B}|) = \frac{2n-1+|B||\mathcal{A}_{\subsetneq B}|+(n-|B|-1)(n-|\mathcal{A}_{\subsetneq B}|-3)}{n}\textrm{.}\]

\medskip
\smallskip

\noindent In the numerator of $\zeta(|B|, |\mathcal{A}_{\subsetneq B}|)$, $2n-1$ comes from equality with $|Y|+n$, $|B||\mathcal{A}_{\subsetneq B}|$ comes from being less than or equal to $\sum_{A \in \mathcal{A}_{\leq |B|}}|A|$ (by Proposition C), and $(n-|B|-1)(n-|\mathcal{A}_{\subsetneq B}|-3)$ comes from all $n-|\mathcal{A}_{\subsetneq B}|-3$ member sets of $\hat{\mathcal{A}}_{> |B|} \setminus \{Y,[n]\}$ having size greater than or equal to $n-|B|-1$. (If there exists $X \in \hat{\mathcal{A}}_{> |B|} \setminus \{Y,[n]\}$ with $|X|<n-|B|-1$, then $|B \cup X| < n-1$ and by Proposition B, we may assume that there exists $A \in \mathcal{A}_{\subsetneq B}$. It follows that $A \subsetneq B \subsetneq B \cup X \subsetneq [n]$, which contradicts Lemma 1.3.) A continuous relaxation of $\zeta$ is the function $f \colon \mathbb{R}^2 \to \mathbb{R}$ such that:

\[\ \ f(x,y) = n-x-y-2 + \frac{2xy+3x+y+2}{n}{.}\] 

\medskip
\smallskip

\noindent For Case 1, $B$ belongs to $\mathcal{A}_{<n/2}$, so $|B| \leq \frac{n-1}{2}$. Also, by Proposition B we may assume that $1 \leq |\mathcal{A}_{\subsetneq B}| \leq |B|$. Therefore, solving the following problem would provide a lower bound for $\textrm{Avg}(\hat{\mathcal{A}})$:

\medskip

\[\ \ \min\Bigr\{f(x,y)\Bigr\} \textrm{ s.t. }1 \leq y \leq x \leq \frac{n-1}{2}\textrm{.}\]

\medskip
\smallskip

\noindent \[\textrm{\hspace{-4.35cm}The problem has solution:} \hspace{0.95cm}f(x^*,y^*)=\frac{n}{2}\textrm{, at }(x^*,y^*) = \ \Bigr(\frac{n}{2}-1, \frac{n}{2}-1\Bigr)\textrm{.}\]

\medskip
\medskip

\noindent Hence, $\textrm{Avg}(\mathcal{A}) \geq \frac{n}{2}$, proving Case 1 of Theorem 2.1. 

\bigskip

\noindent Denote by $\binom{S}{k}$ the family of $k$-element subsets of a set $S$. We conjecture that the minimum average member set size of a union-closed family $\mathcal{A}$ such that $h = 4 \leq n$ and $|\mathcal{B}|=1$ is achieved by the following family $\mathcal{A}^*$ with $b(\mathcal{A}^*) = [n]$, as illustrated in Figure 2.1:

\[\mathcal{A}^{*} = \ \Bigr\{[n], \Bigr[\Bigr\lceil\frac{n}{2}\Bigr\rceil-1\Bigr] \Bigr\} \ \ \cup \ \ \Bigr\{[n]\setminus\{x\} \ \Bigr| \ x \in \ [n] \ \setminus \ \Bigr[\Bigr\lceil\frac{n}{2}\Bigr\rceil \Bigr]\Bigr\} \ \ \cup \ \ \binom{[\lceil\frac{n}{2}\rceil - 1]}{\lceil\frac{n}{2}\rceil-2}\textrm{.}\]

\medskip

\noindent To verify that $\mathcal{A}^*$ is separating, we observe that for any two distinct elements $x$ and $y$ in $[n]$, the following member set $S$ of $\mathcal{A}^*$ contains exactly one element from $\{x,y\}$:

\smallskip
\smallskip

\[S=\begin{cases}

[\lceil\frac{n}{2}\rceil-1] \setminus \{x\} & \textrm{ if } x,y \in [\lceil\frac{n}{2}\rceil-1] \\

[\lceil\frac{n}{2}\rceil-1] & \textrm{ if } \min\{x,y\} \in [\lceil\frac{n}{2}\rceil-1] \ \land \ \max\{x,y\} \in [n] \setminus [\lceil\frac{n}{2}\rceil-1] \\

[n]\setminus\{\max\{x,y\}\} & \textrm{ if } x, y \in [n] \setminus [\lceil\frac{n}{2}\rceil-1]

\end{cases}\textrm{.}\]

\begin{figure}

\caption*{\textbf{Figure 2.1:} Both parameters of $\zeta$ are equal to $\lceil\frac{n}{2}\rceil-1$ for $\mathcal{A}^*$.}

\vspace{-0.25cm}

\begin{center}

\begin{tikzpicture}[scale=0.84]

\node[scale=0.8375] (a) at (-5.75, 2.375) {};

\node[scale=0.8375] (b) at (-5.75, -9.875) {};

\draw (a)--(b);

\node[scale=0.8375] (c) at (-8.1, 1.25) {};

\node[scale=0.8375] (d) at (5.25, 1.25) {};

\draw (c)--(d);

\node[scale=0.8375] at (-0.25,1.725) {\large{$A \in \mathcal{A}^*$}};

\node[scale=0.8375] at (-6.9,1.75) {\large{$|A|$}};

\node[scale=0.8375] at (-6.96,0.44) {\large{$n$}};

\node[scale=0.8375] at (-7,-0.078) {\large{$n-1$}};

\node[scale=0.8375] at (-7,-0.57) {\large{$n-1$}};

\node[scale=0.8375] at (-7,-1.04) {\large{$n-1$}};

\node[scale=0.8375,
    minimum width = 1cm, 
    minimum height = 0.5cm] at (-6.95,-1.8225) {\Huge{$\vdots$}};

\node[scale=0.8375] at (-7,-3.05) {\large{$n-1$}};

\node[scale=0.8375] at (-7,-3.545) {\large{$n-1$}};

\node[scale=0.8375] at (-7,-4.0475) {\large{$n-1$}};

\node[scale=0.8375] at (-7.1,-4.5475) {\large{$\lceil n/2 \rceil-1$}};

\node[scale=0.8375] at (-7.1,-5.0475) {\large{$\lceil n/2 \rceil-2$}};

\node[scale=0.8375] at (-7.1,-5.575) {\large{$\lceil n/2 \rceil-2$}};

\node[scale=0.8375] at (-7.1,-6.075) {\large{$\lceil n/2 \rceil-2$}};

\node[scale=0.8375,
    minimum width = 0.5cm,
    minimum height = 0.5cm] at (-6.95,-6.85) {\Huge{$\vdots$}};

\node[scale=0.8375] at (-7.1,-8.15) {\large{$\lceil n/2 \rceil-2$}};

\node[scale=0.8375] at (-7.1,-8.7) {\large{$\lceil n/2 \rceil-2$}};

\node[scale=0.8375] at (-7.1,-9.25) {\large{$\lceil n/2 \rceil-2$}};

\node[scale=0.8375,rectangle,
    draw = black,
    text = black,
    fill = lightgray,
    minimum width = 10cm, 
    minimum height = 0.5cm] at (-0.275,0.5) {};

\node[scale=0.8375,rectangle,
    draw = black,
    text = black,
    fill = lightgray,
    minimum width = 9.5cm, 
    minimum height = 0.5cm] at (-0.525,0) {};

\node[scale=0.8375, rectangle,
    draw = black,
    text = black,
    fill = lightgray,
    minimum width = 9cm, 
    minimum height = 0.5cm] at (-0.775,-0.5) {};

\node[scale=0.8375, rectangle,
    draw = black,
    text = black,
    fill = lightgray,
    minimum width = 0.5cm, 
    minimum height = 0.5cm] at (4.475,-0.5) {};
      
\node[scale=0.8375, rectangle,
    draw = black,
    text = black,
    fill = lightgray,
    minimum width = 8.5cm,
    minimum height = 0.5cm] at (-1.025,-1) {};

\node[scale=0.8375, rectangle,
    draw = black,
    text = black,
    fill = lightgray,
    minimum width = 1cm, 
    minimum height = 0.5cm] at (4.225,-1) {};  

\node[scale=0.8375,
    minimum width = 1cm, 
    minimum height = 0.5cm] at (-0.275,-1.775) {\Huge{$\vdots$}};

\node[scale=0.8375, rectangle,
    draw = black,
    text = black,
    fill = lightgray,
    minimum width = 6cm,
    minimum height = 0.5cm] at (-2.275,-3) {};

\node[scale=0.8375, rectangle,
    draw = black,
    text = black,
    fill = lightgray,
    minimum width = 3.5cm,
    minimum height = 0.5cm] at (2.975,-3) {};

\node[scale=0.8375, rectangle,
    draw = black,
    text = black,
    fill = lightgray,
    minimum width = 5.5cm,
    minimum height = 0.5cm] at (-2.525,-3.5) {};

\node[scale=0.8375, rectangle,
    draw = black,
    text = black,
    fill = lightgray,
    minimum width = 4cm,
    minimum height = 0.5cm] at (2.725,-3.5) {}; 

\node[scale=0.8375, rectangle,
    draw = black,
    text = black,
    fill = lightgray,
    minimum width = 5cm, 
    minimum height = 0.5cm] at (-2.775,-4) {};

\node[scale=0.8375, rectangle,
    draw = black,
    text = black,
    fill = lightgray,
    minimum width = 4.5cm,
    minimum height = 0.5cm] at (2.475,-4) {}; 

\node[scale=0.8375, rectangle,
    draw = black,
    text = black,
    fill = lightgray,
    minimum width = 4.5cm, 
    minimum height = 0.5cm] at (-3.025,-4.5) {};

\node[scale=0.8375, rectangle,
    draw = black,
    text = black,
    fill = lightgray,
    minimum width = 4cm, 
    minimum height = 0.5cm] at (-3.275,-5) {};

\node[scale=0.8375, rectangle,
    draw = black,
    text = black,
    fill = lightgray,
    minimum width = 3.5cm, 
    minimum height = 0.5cm] at (-3.525,-5.5) {};

\node[scale=0.8375, rectangle,
    draw = black,
    text = black,
    fill = lightgray,
    minimum width = 0.5cm, 
    minimum height = 0.5cm] at (-1.025,-5.5) {};

\node[scale=0.8375, rectangle,
    draw = black,
    text = black,
    fill = lightgray,
    minimum width = 3cm, 
    minimum height = 0.5cm] at (-3.775,-6) {};

\node[scale=0.8375, rectangle,
    draw = black,
    text = black,
    fill = lightgray,
    minimum width = 1cm, 
    minimum height = 0.5cm] at (-1.275,-6) {};

\node[scale=0.8375,
    minimum width = 0.5cm, 
    minimum height = 0.5cm] at (-2.95,-6.85) {\Huge{$\vdots$}};

\node[scale=0.8375, rectangle,
    draw = black,
    text = black,
    fill = lightgray,
    minimum width = 1cm, 
    minimum height = 0.5cm] at (-4.775,-8.15) {};

\node[scale=0.8375, rectangle,
    draw = black,
    text = black,
    fill = lightgray,
    minimum width = 3cm, 
    minimum height = 0.5cm] at (-2.275,-8.15) {};

\node[scale=0.8375, rectangle,
    draw = black,
    text = black,
    fill = lightgray,
    minimum width = 0.5cm, 
    minimum height = 0.5cm] at (-5.025,-8.65) {};

\node[scale=0.8375, rectangle,
    draw = black,
    text = black,
    fill = lightgray,
    minimum width = 3.5cm, 
    minimum height = 0.5cm] at (-2.525,-8.65) {};

\node[scale=0.8375, rectangle,
    draw = black,
    text = black,
    fill = lightgray,
    minimum width = 4cm, 
    minimum height = 0.5cm] at (-2.775,-9.15) {};

\end{tikzpicture}

\end{center}

\end{figure}
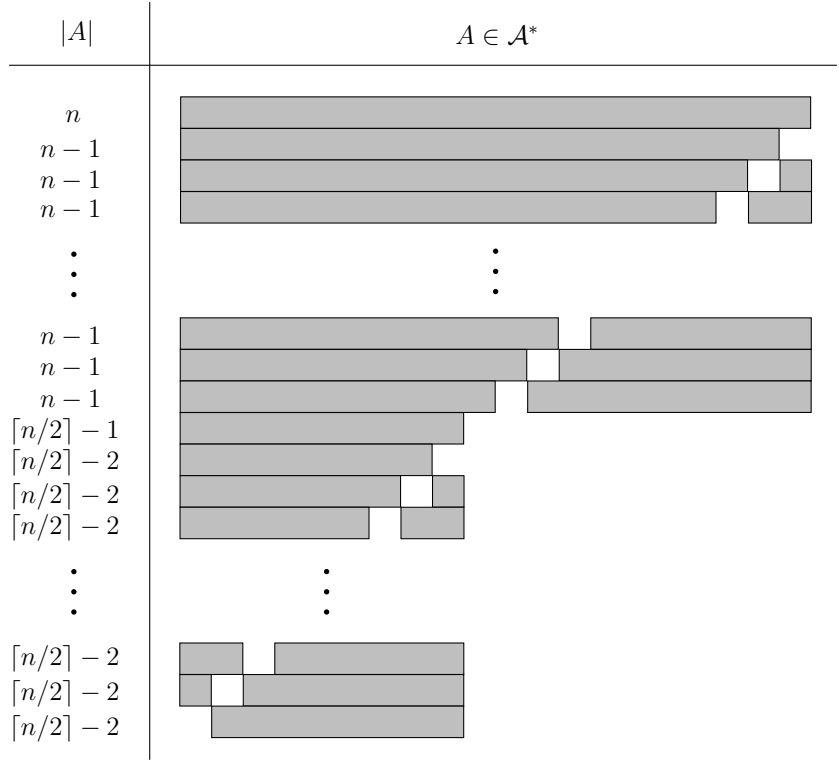

\bigskip
\smallskip

\noindent \textbf{\normalsize{Case 2 $(|\mathcal{B}|=2)$}}
\medskip

\noindent In this case, we let $\mathcal{B}=\{B_1,B_2\}$. We first assume that $|B|<n-1$. 

\medskip

\noindent Let $\hat{\mathcal{A}}=\mathcal{A}_{\subseteq B} \cup \{Y,[n]\}$. We observe that $|B| \geq \frac{n}{2}$, and that $|A| \geq \frac{n}{2}$ for any $A \in \mathcal{A}_{< |B|} \setminus \mathcal{A}_{\subsetneq B}$. Therefore, any member set of $\mathcal{A} \setminus \hat{\mathcal{A}}$ has size greater than or equal to $\frac{n}{2}$, and it is sufficient for resolving this subcase to prove that $\textrm{Avg}(\hat{\mathcal{A}}) \geq \frac{n}{2}$. We introduce the lower bound $\eta(|B|, |\mathcal{A}_{\subsetneq B}|)$ of $\textrm{Avg}(\hat{\mathcal{A}})$:

\medskip

\[\hspace{-0.875cm} \textrm{Avg}(\hat{\mathcal{A}})=\frac{\sum_{A \in \hat{\mathcal{A}}}|A|}{|\hat{\mathcal{A}}|} = \frac{|Y|+n+\sum_{A \in \mathcal{A}_{\subseteq B}}|A|}{|\mathcal{A}_{\subsetneq B}|+3}\]

\smallskip

\[\hspace{0.185cm} \geq \eta(|B|, |\mathcal{A}_{\subsetneq B}|) = \frac{2n-1+|B||\mathcal{A}_{\subsetneq B}|}{|\mathcal{A}_{\subsetneq B}|+3}\textrm{.}\]

\medskip
\smallskip

\noindent In the numerator of $\eta(|B|, |\mathcal{A}_{\subsetneq B}|)$, $2n-1$ again comes from being equal to $|Y|+n$, and $|B||\mathcal{A}_{\subsetneq B}|$ comes from being less than or equal to $\sum_{A \in \mathcal{A}_{\subseteq B}}|A|$. Then similar to $f$ from Case 1, a continuous relaxation of $\eta$ is the function $g \colon \mathbb{R} \times \mathbb{R} \setminus \{-3\}  \to \mathbb{R}$ such that:

\medskip

\[\ \ g(x,y) = \frac{2n+xy-1}{y+3}{.}\] 

\medskip
\smallskip

\noindent Recall that $\frac{n}{2} \leq |B| < n-1$, and by Proposition B, we may assume that $1 \leq |\mathcal{A}_{\subsetneq B}| \leq |B|$. Therefore, solving the following problem would provide a lower bound for $\textrm{Avg}(\hat{\mathcal{A}})$:

\medskip

\[\ \ \min\Bigr\{g(x,y)\Bigr\} \textrm{ s.t. }\Bigr(\frac{n}{2} \leq x \leq n-2\Bigr) \ \land \ \Bigr(1 \leq y \leq x\Bigr)\textrm{.}\]

\medskip
\smallskip

\noindent \[\textrm{\hspace{-3.95cm}The problem has solution:}\hspace{0.525cm}g(x^*,y^*)=\frac{n}{2}+\frac{n-2}{n+6}>\frac{n}{2}\textrm{, at }(x^*,y^*) = \ \Bigr(\frac{n}{2}, \frac{n}{2}\Bigr)\textrm{.}\]

\medskip
\smallskip

\noindent Thus, $\textrm{Avg}(\mathcal{A}) \geq \frac{n}{2}$, resolving Case 2 for $|B| < n-1$.

\medskip

\noindent Next, we assume that $|B|=n-1$, and let $B_1=[\frac{n-1}{2}]$ and $B_2=[n-1] \setminus B_1$ without loss of generality. 

\medskip

\noindent If $\mathcal{A}_{<n/2}=\mathcal{B}$, then we let $\hat{\mathcal{A}}=\{B_1,B_2,[n]\}$, so $\textrm{Avg}(\hat{\mathcal{A}})=\frac{2n-1}{3} \geq \frac{n}{2}$. Then $\textrm{Avg}(\mathcal{A}) \geq \frac{n}{2}$ because any member set in $\mathcal{A} \setminus \hat{\mathcal{A}}$ has size greater than or equal to $\frac{n}{2}$.

\medskip

\noindent If there exist $X_1$ and $X_2$ in $\mathcal{A}_{<n/2}$ such that $(X_1 \subsetneq B_1) \land (X_2 \subsetneq B_2)$, then $X_1 \subsetneq B_1 \subsetneq B_1 \cup X_2 \subsetneq B \subsetneq [n]$, forming a chain of size $5$ in $\mathcal{A}$, which contradicts $h=4$.

\pagebreak

\noindent There are two remaining assumptions that are possible for Case 2 when $|B|=n-1$:

\bigskip

\noindent 1.) $\mathcal{A}_{<n/2} \neq \mathcal{B} \ \land \ ((A \in \mathcal{A}_{<n/2} \setminus \mathcal{B} \implies A \subsetneq B_1) \lor (A \in \mathcal{A}_{<n/2} \setminus \mathcal{B} \implies A \subsetneq B_2))$:

\medskip
\smallskip

\noindent Without loss of generality, we assume that $A \in \mathcal{A}_{<n/2} \setminus \mathcal{B} \implies A \subsetneq B_1$. We observe that $\mathcal{A}'=\mathcal{A} \setminus \{B_2\}$ is a union-closed family with $\mathcal{B}(\mathcal{A}')=\{B_1\}=\{B'\}$. Additionally, we have that $\mathcal{A}'$ is separating. (Otherwise, there exist $x \in B_2$ and $y \in [n] \setminus B_2$ such that for any $A \in \mathcal{A}'$, $(x \in A \ \lor \ y \in A) \implies (x \in A \ \land \ y \in A)$. Then $B_1 \cap B_2 = \emptyset$ implies that $y \not \in B_1$, making $y=n$. Thus, $y \not \in [n-1]$, yet $x \in [n-1] = B_1 \cup B_2 \in \mathcal{A}'$, a contradiction.) Hence, $|\mathcal{A}'| \geq n$ by Lemma 2.1.1, and $\mathcal{A}'$ falls within the scope of Case 1. We apply the proof for Case 1 to $\mathcal{A}'$, considering the corresponding family $\hat{\mathcal{A}}'=\mathcal{A}' \setminus \tilde{\mathcal{A}}'$, where $\tilde{\mathcal{A}}' \subseteq \mathcal{A}'_{>|B'|} \setminus \{Y, [n]\}$ such that $|\tilde{\mathcal{A}}'|=|\mathcal{A}'|-n$. Because $|B_2|=\frac{n-1}{2}=n-|B'|-1$, $B_2$ may be counted as one of the $n-|{\mathcal{A}'}_{\subsetneq B'}|-3$ member sets of $\hat{\mathcal{A}}'_{>|B'|} \setminus \{Y,[n]\}$ of size greater than or equal to $n-|B'|-1$. (If there were no such member set, then $n-|\mathcal{A}'_{\subsetneq B'}|-3 = 0 \geq n-\frac{n-1}{2}-3$, implying that $n \leq 5$. Then $n$ being odd would imply that $n=5$, and $\mathcal{A}'$ would be equal to $\{\{1,2,3,4,5\}, Y, \{1,2\}, \{1\}, \{2\}\}$ for some $Y$ of size $4$, contradicting the fact that $\mathcal{A}'$ is separating.) It is thus sufficient for resolving this subcase to show that $f(\frac{n-1}{2},y) \geq \frac{n}{2}$ for $1 \leq y \leq \frac{n-1}{2}$, which follows from $f(\frac{n-1}{2},y)$ being greater than or equal to $f(\frac{n}{2}-1,\frac{n}{2}-1)$ for all such $y$.

\bigskip

\noindent 2.) There exists $B_3$ in $\mathcal{A}_{<n/2}$ such that $(B_1 \cap B_3 \neq \emptyset) \land (B_2 \cap B_3 \neq \emptyset)$:

\medskip
\smallskip

\noindent If $\mathcal{A}_{<n/2}=\{B_1,B_2,B_3\}$, then $\textrm{Avg}(\hat{\mathcal{A}}) \geq \frac{2n+1}{4}$ for $\hat{\mathcal{A}}=\{B_1,B_2,B_3,[n]\}$, making $\textrm{Avg}(\mathcal{A}) \geq \frac{n}{2}$, as any member set of $\mathcal{A} \setminus \hat{\mathcal{A}}$ must have size greater than or equal to $\frac{n}{2}$. Else, $\{B_1,B_2,B_3\} \subsetneq \mathcal{A}_{<n/2}$, and we consider the partition $\mathcal{P}=\{P_1,P_2,P_3,P_4\}$ of $B$, where $P_1 = B_1 \setminus B_3$, $P_2 = B_3 \setminus B_2$, $P_3 = B_3 \setminus B_1$, and  $P_4 = B_2 \setminus B_3$. For any $X \in \mathcal{A}_{<n/2} \setminus \{B_1,B_2,B_3\}$ and $i \in \{1,2,3,4\}$, we have that $X \cap P_i \in \{\emptyset, P_i\}$. (If not, then there exists $X \in \mathcal{A}_{<n/2} \setminus \{B_1,B_2,B_3\}$ such that, without loss of generality, $\emptyset \subsetneq X \cap P_1 \subsetneq P_1$, and we have that $B_2 \subsetneq B_2 \cup B_3 \subsetneq B_2 \cup B_3 \cup X \subsetneq B \subsetneq [n]$, which forms a chain in $\mathcal{A}$ of size $5$, contradicting $h=4$.) Now, we observe that $|\bigcup_{P \in \mathcal{P} \setminus \{P_i\}}P| \geq \frac{n+1}{2}$ for any $i \in \{1,2,3,4\}$. We also observe that $(\mathcal{P} \cup \{\emptyset\}) \cap (\mathcal{A}_{<n/2} \setminus \{B_1,B_2,B_3\}) = \emptyset$. (Otherwise, there exists $X \in \mathcal{A}_{<n/2} \setminus \{B_1,B_2,B_3\}$ such that either $X=\emptyset \subsetneq B_1 \subsetneq B_1 \cup B_3 \subsetneq B \subsetneq [n]$, or without loss of generality, $X=P_1 \subsetneq B_1 \subsetneq B_1 \cup B_3 \subsetneq B \subsetneq [n]$, again contradicting $h=4$.) It follows that any $X \in \mathcal{A}_{<n/2}$ must be the union of exactly two distinct member sets of $\mathcal{P}$. The maximum cardinality of $\mathcal{A}_{<n/2}$ is thus $\binom{|\mathcal{P}|}{2}=6$, as illustrated in Figure 2.2.

\medskip

\begin{figure}[H]

\caption*{\textbf{Figure 2.2:} $\mathcal{A}_{<n/2} = \{B_1,B_2,B_3,B_4,B_5,B_6\}$ has maximum size.}

\vspace{-0.25cm}

\begin{center}

\begin{tikzpicture}

\draw [decorate,decoration={brace,amplitude=7pt},xshift=-4pt,yshift=0pt]
(-2.5,0.65) -- (7.725,0.65) node [black,midway,xshift=-0.6cm]{};

\node at (2.275,1.275) {\large{$B=[n-1]$}};

\node at (-3.15,0) {\large{$B_1 \ \Bigr\{$}};
\node at (-3.15,-0.75) {\large{$B_2 \ \Bigr\{$}};
\node at (-3.15,-1.5) {\large{$B_3 \ \Bigr\{$}};
\node at (-3.15,-2.25) {\large{$B_4 \ \Bigr\{$}};
\node at (-3.15,-3) {\large{$B_5 \ \Bigr\{$}};
\node at (-3.15,-3.75) {\large{$B_6 \ \Bigr\{$}};

\node[rectangle,
    draw = black,
    text = black,
    fill = lightgray,
    minimum width = 5cm, 
    minimum height = 0.75cm] at (0,0) {};

\node[rectangle,
    draw = black,
    text = black,
    fill = lightgray,
    minimum width = 5cm, 
    minimum height = 0.75cm] at (5,-0.75) {};

\node[rectangle,
    draw = black,
    text = black,
    fill = lightgray,
    minimum width = 5cm, 
    minimum height = 0.75cm] at (2.5,-1.5) {};

\node[rectangle,
    draw = black,
    text = black,
    fill = lightgray,
    minimum width = 2.5cm, 
    minimum height = 0.75cm] at (-1.25,-2.25) {};
    
\node[rectangle,
    draw = black,
    text = black,
    fill = lightgray,
    minimum width = 2.5cm, 
    minimum height = 0.75cm] at (6.25,-2.25) {};

\node[rectangle,
    draw = black,
    text = black,
    fill = lightgray,
    minimum width = 2.5cm, 
    minimum height = 0.75cm] at (1.25,-3) {};
    
\node[rectangle,
    draw = black,
    text = black,
    fill = lightgray,
    minimum width = 2.5cm, 
    minimum height = 0.75cm] at (6.25,-3) {};

\node[rectangle,
    draw = black,
    text = black,
    fill = lightgray,
    minimum width = 2.5cm, 
    minimum height = 0.75cm] at (-1.25,-3.75) {};
    
\node[rectangle,
    draw = black,
    text = black,
    fill = lightgray,
    minimum width = 2.5cm, 
    minimum height = 0.75cm] at (3.75,-3.75) {};

\draw [rotate=180,decorate,decoration={brace,amplitude=7pt},xshift=-4pt,yshift=0pt]
(0.075,4.425) -- (2.675,4.425) node [black,midway,xshift=-0.6cm]{};

\node at (-1.225,-5.075) {\large{$P_1$}};

\draw [rotate=180,decorate,decoration={brace,amplitude=7pt},xshift=-4pt,yshift=0pt]
(-2.375,4.425) -- (0.075,4.425) node [black,midway,xshift=-0.6cm]{};

\node at (1.3,-5.075) {\large{$P_2$}};

\draw [rotate=180,decorate,decoration={brace,amplitude=7pt},xshift=-4pt,yshift=0pt]
(-4.9,4.425) -- (-2.375,4.425) node [black,midway,xshift=-0.6cm]{};

\node at (3.8,-5.075) {\large{$P_3$}};

\draw [rotate=180,decorate,decoration={brace,amplitude=7pt},xshift=-4pt,yshift=0pt]
(-7.425,4.425) -- (-4.9,4.425) node [black,midway,xshift=-0.6cm]{};

\node at (6.3,-5.075) {\large{$P_4$}};

\end{tikzpicture}

\end{center}

\end{figure}
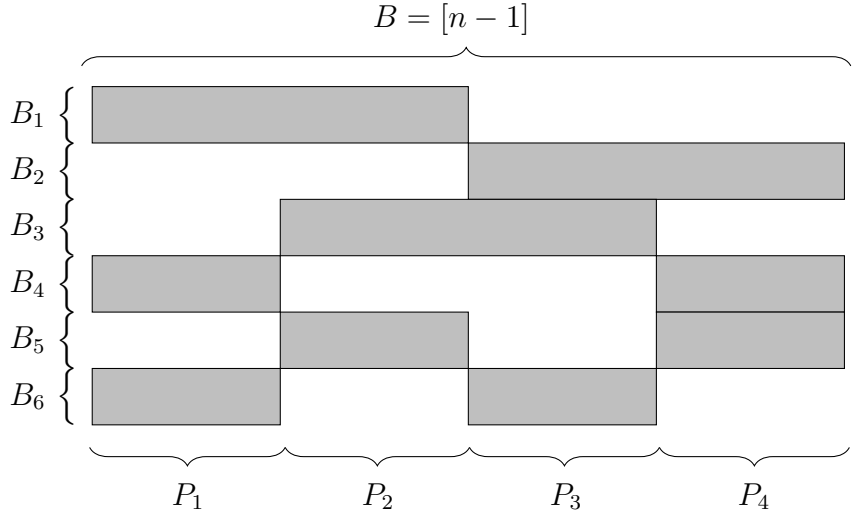

\smallskip

\noindent Therefore, we have that $|\mathcal{A}_{<n/2}| \in \{4,5,6\}$. We next establish Proposition D and, under this present assumption, Proposition E.

\pagebreak

\noindent \textbf{Proposition D.} If $N$ is a positive integer, then for any $P = \{p_1, \cdots, p_N\} \in \mathbb{R}^N$ and $k \in [N]$:

\[\binom{N-1}{k-1} \sum_{i \in [N]} p_i \ = \sum_{S \in \binom{[N]}{k}} \ \sum_{j \in S} p_j\textrm{.}\]

\noindent \textit{Proof.} For any $m \in [N]$, there are $\binom{N-1}{k-1}$ subsets of $P$ that contain both $p_m$ and exactly $k-1$ other elements of $P$, making $p_m$ occur on the right-hand side $\binom{N-1}{k-1}$ times. The proposition then follows from adding these occurences together across all $m \in [N]$.

\medskip

\noindent We observe that setting $P$ from Proposition D equal to $\{1\}^N$ yields the identity $\binom{N-1}{k-1}N = \binom{N}{k} k$.

\medskip

\noindent \textbf{Proposition E.} For any four distinct $A_1$, $A_2$, $A_3$, and $A_4$ in $\mathcal{A}_{<n/2}$, $\sum_{i=1}^4 |A_i| \geq \frac{3n+1}{2}$.

\medskip

\noindent \textit{Proof.} There must exist two member sets from $\{A_1, A_2, A_3, A_4\}$ whose union is equal to $[n-1]$. Thus, we assume without loss of generality that $A_1 = B_1$, $A_2=B_2$, and $A_3 = B_3 = P_2 \cup P_3$. If $A_4=P_1 \cup P_4$, then $|A_1|+|A_2|+|A_3|+|A_4|=2n-2 \geq \frac{3n+1}{2}$. Else, we have that $A_4=P_1 \cup P_3$ and $B_1 \subsetneq A_3 \cup A_4$, or otherwise that $A_4=P_2 \cup P_4$ and $B_2 \subsetneq A_3 \cup A_4$. In either case, $A_3 \cap A_4 \neq \emptyset$, and $(|A_1|+|A_2|)+(|A_3|+|A_4|) \geq (n-1)+(\frac{n-1}{2}+2) = \frac{3n+1}{2}$, completing the proof of Proposition E.

\medskip

\noindent We let $\hat{\mathcal{A}}=\mathcal{A}_{<n/2} \cup \{Y,[n]\}$ and address the three possible values of $|\mathcal{A}_{<n/2}|$:

\bigskip

\noindent (\romannumeral 1) $|\mathcal{A}_{<n/2}|=4$: Let $\mathcal{A}_{<n/2}=\{B_1,B_2,B_3,B_4\}$. By Proposition E, we have that: 

\[\textrm{Avg}(\hat{\mathcal{A}}) = \frac{|Y|+n+\sum_{i=1}^4|B_i|}{6} \geq \frac{2n-1+\frac{3n+1}{2}}{6} = \frac{7n-1}{12} > \frac{n}{2}\textrm{.}\]

\bigskip
\medskip

\noindent (\romannumeral 2) $|\mathcal{A}_{<n/2}|=5$: Let $\mathcal{A}_{<n/2}=\{B_1,B_2,B_3,B_4,B_5\}$. Applying Proposition D, we set $N=5$, $k=4$, and $p_i=|B_i|$ for each $i \in \{1,2,3,4,5\}$ to obtain:

\[\binom{4}{3} \sum_{i=1}^5 |B_i| \ = \ \sum_{S \in \binom{[5]}{4}} \ \sum_{j \in S} |B_j|\textrm{.}\]

\noindent 

\[\hspace{-5.825cm} \textrm{It follows by Proposition E that: }\hspace{0.585cm} \binom{4}{3} \sum_{i=1}^5 |B_i| \ \geq \ \binom{5}{4} \frac{3n+1}{2}\textrm{.}\]

\medskip
\smallskip

\[\hspace{-2.825cm} \textrm{Thus, we have that:} \hspace{0.35cm} \textrm{Avg}(\hat{\mathcal{A}}) = \frac{|Y|+n+\sum_{i=1}^5 |B_i|}{7} \geq \frac{2n-1+\frac{15n+5}{8}}{7} = \frac{31n-3}{56} > \frac{n}{2}\textrm{.}\]

\bigskip
\medskip

\noindent (\romannumeral 3) $|\mathcal{A}_{<n/2}|=6$: Let $\mathcal{A}_{<n/2}=\{B_1,B_2,B_3,B_4,B_5,B_6\}$. We again apply Proposition D, this time setting $N=6$, $k=4$, and $p_i=|B_i|$ for each $i \in \{1,2,3,4,5,6\}$ to obtain:

\[\binom{5}{3} \sum_{i = 1}^6 |B_i| \ = \ \sum_{S \in \binom{[6]}{4}} \ \sum_{j \in S} |B_j|\textrm{.}\]

\medskip

\[\hspace{-5.75cm} \textrm{By Proposition E, we have that: }\hspace{0.65cm} \binom{5}{3} \sum_{i=1}^6 |B_i| \ \geq \ \binom{6}{4} \frac{3n+1}{2}\textrm{.}\]

\medskip

\[\hspace{-3.25cm} \textrm{As a result:} \hspace{1.125cm} \textrm{Avg}(\hat{\mathcal{A}}) = \frac{|Y|+n+\sum_{i=1}^6 |B_i|}{8} \geq \frac{2n-1+\frac{9n+3}{4}}{8} = \frac{17n-1}{32} > \frac{n}{2}\textrm{.}\]

\medskip

\noindent Thus, $\textrm{Avg}(\hat{\mathcal{A}}) \geq \frac{n}{2}$ for any $|\mathcal{A}_{<n/2}|$. It follows that $\textrm{Avg}(\mathcal{A}) \geq \frac{n}{2}$, as any member set of $\mathcal{A} \setminus \hat{\mathcal{A}}$ has size greater than or equal to $\frac{n}{2}$, proving the subcase $|B|=n-1$ of Case 2.

\bigskip

\noindent This concludes the proof of Case 2, which then completes the overall proof of Theorem 2.1.

\bigskip
\medskip

\noindent \textbf{Corollary 2.2.} \textit{For any separating union-closed family} $\mathcal{A}$ \textit{with} $h = 4 \leq n$ \textit{and} $0 \leq |\mathcal{B}| \leq 2$\textit{, there exists an element of} $[n]$ \textit{that is in at least} $\frac{|\mathcal{A}|}{2}$ \textit{member sets of} $\mathcal{A}$\textit{.}

\medskip

\noindent \textit{Proof.} If not, then there is a separating union-closed family $\mathcal{A}$ with $h = 4 \leq n$ and $0 \leq |\mathcal{B}| \leq 2$ such that the family $\mathcal{A}_{\{x\}}=\{A \in \mathcal{A} \ | \ x \in A\}$ has size less than $\frac{|\mathcal{A}|}{2}$ for every $x \in [n]$. It follows that $\textrm{Avg}(\mathcal{A}) = \frac{\sum_{A \in \mathcal{A}}|A|}{|\mathcal{A}|} = \frac{\sum_{x \in [n]}|\mathcal{A}_{\{x\}}|}{|\mathcal{A}|} < \frac{n(|\mathcal{A}| / 2)}{|\mathcal{A}|} = \frac{n}{2}$. This contradicts Theorem 2.1, which states that $\textrm{Avg}(\mathcal{A}) \geq \frac{n}{2}$.

\section*{3. Limitation of averaging for larger values of $h$}

In this section, Theorem 3.2 demonstrates that the technique of averaging member set sizes of separating union-closed families cannot be used to prove Conjecture 1.1 for any $h \geq 5$. We first establish Lemma 3.1.

\medskip
\smallskip

\noindent \textbf{Lemma 3.1.} \textit{For any integer} $n \geq 9$\textit{, there exists a separating union-closed family} $\mathcal{A}$ \textit{with} $h=5$ \textit{and} $|\mathcal{B}|=1$ \textit{such that}:

\vspace{-0.25cm}

\[\textrm{Avg}(\mathcal{A}) = \frac{\sum_{A \in \mathcal{A}}|A|}{|\mathcal{A}|} < \frac{n}{2}\textrm{.}\]

\medskip
\medskip

\noindent \textit{Proof.} Consider $\mathcal{A}^*$ from Figure 2.1. For any integer $n \geq 9$, we construct a union-closed family $\mathcal{A}^{**}$ with base set $[n]$, height $h=5$, and $|\mathcal{B}(\mathcal{A}^{**})|=1$ as follows:

\smallskip

\[\mathcal{A}^{**} = \mathcal{A}^{*} \cup \binom{[\lceil \frac{n}{2} \rceil-1]}{\lceil \frac{n}{2} \rceil-3} = \ \Bigr\{[n], \Bigr[\Bigr\lceil \frac{n}{2} \Bigr\rceil-1\Bigr] \Bigr\} \ \ \cup \ \ \Bigr\{[n]\setminus\{x\} \ \Bigr| \ x \in \ [n] \ \setminus \ \Bigr[\Bigr\lceil \frac{n}{2} \Bigr\rceil\Bigr]\Bigr\} \ \ \cup \ \ \bigcup_{i=2}^3\binom{[\lceil \frac{n}{2} \rceil - 1]}{\lceil \frac{n}{2} \rceil-i}\textrm{.}\]

\medskip
\smallskip

\noindent $\mathcal{A}^{**}$ is separating because it is a superfamily of $\mathcal{A}^*$, which itself is a separating union-closed family with $b(\mathcal{A}^*)=b(\mathcal{A}^{**})=[n]$. We compute the average size of a member set from $\mathcal{A}^{**}$ to be:

\bigskip
\smallskip

\noindent \[\textrm{Avg}(\mathcal{A}^{**}) = \frac{n+\Bigr(\lceil \frac{n}{2} \rceil-1\Bigr)+\Bigr(n - \lceil\frac{n}{2} \rceil\Bigr)\Bigr(n-1\Bigr)+\Bigr(\lceil \frac{n}{2} \rceil-1\Bigr)\Bigr(\lceil \frac{n}{2} \rceil-2\Bigr)+\frac{1}{2}\Bigr(\lceil \frac{n}{2} \rceil-1\Bigr)\Bigr(\lceil \frac{n}{2} \rceil-2\Bigr)\Bigr(\lceil \frac{n}{2} \rceil-3\Bigr)}{2+\Bigr(n - \lceil \frac{n}{2} \rceil\Bigr)+\Bigr(\lceil \frac{n}{2} \rceil-1\Bigr)+\frac{1}{2}\Bigr(\lceil \frac{n}{2} \rceil-1\Bigr)\Bigr(\lceil \frac{n}{2} \rceil-2\Bigr)}\textrm{.}\]

\medskip
\smallskip

\noindent Thus, $\textrm{Avg}(\mathcal{A}^{**})$ is equal to $\frac{n^3+36n-32}{2n^2+4n+32}$ when $n$ is even, and equal to $\frac{n^3+3n^2+15n-3}{2n^2+8n+22}$ when $n$ is odd. In either case, $\textrm{Avg}(\mathcal{A}^{**}) < \frac{n}{2}$, completing the proof of Lemma 3.1.

\bigskip

For proving Theorem 3.2, we let $[0]=\emptyset$ and $[c,d]=\{i \in \mathbb{Z} \ | \ c \leq i \leq d\}$ for any two integers $c$ and $d$.

\bigskip

\noindent \textbf{Theorem 3.2.} \textit{For any integer} $n \geq 11$ \textit{and} $k \in [5,n+1]$\textit{, there exists a separating union-closed family} $\mathcal{A}^{(k)}$ \textit{with base set} $[n]$ \textit{and height} $h=k$\textit{, as well as} $|\mathcal{B}(\mathcal{A}^{(k)})|=1$\textit{, such that}:

\[\textrm{Avg}(\mathcal{A}^{(k)}) = \frac{\sum_{A \in \mathcal{A}^{(k)}}|A|}{|\mathcal{A}^{(k)}|} < \frac{n}{2}\textrm{.}\]

\medskip

\noindent \textit{Proof.} For any integer $n \geq 11$, we set $\mathcal{A}^{(5)}=\mathcal{A}^{**}$ and $\Delta = n - 2\lceil \frac{n}{2} \rceil +2$, and establish a conditional recurrence relation for $k \in [5, n]$ as follows:

\smallskip
 
\[\mathcal{A}^{(k+1)} = \mathcal{A}^{(k)} \ \cup \ \begin{cases}  

\Bigr\{\Bigr[\lceil \frac{n}{2} \rceil+k-5 \Bigr]\Bigr\} & \textrm{if }k \in \ \Bigr[5, 5+\Delta\Bigr] \\[5pt]
\Bigr\{\Bigr[\lceil \frac{n}{2} \rceil-\frac{k+2-\Delta}{2}\Bigr]\Bigr\} & \textrm{if } k \in \ \Bigr\{2i+\Delta \ | \ i \in \ \Bigr[3, \frac{n-\Delta}{2}\Bigr]\Bigr\} \\[5pt] 
\Bigr\{\Bigr[\lceil \frac{n}{2} \rceil + \frac{k-5+\Delta}{2}\Bigr]\Bigr\} & \textrm{if }k \in \ \Bigr\{2i+\Delta+1 \ | \ i \in \ \Bigr[3, \frac{n-\Delta-2}{2}\Bigr]\Bigr\}\end{cases}\textrm{.}\]

\medskip

\noindent In the same way that $\mathcal{A}^{**}$ was separating in the proof of Lemma 3.1, $\mathcal{A}^{(k)}$ is separating because it is a superfamily of $\mathcal{A}^*$, which itself is a separating family with the same base set.

\medskip

\noindent For $k \in [6, 6+\Delta]$, we have that $\mathcal{A}^{(k)}=\mathcal{A}^{**} \cup \ \bigcup_{i=5}^{k-1} \{[\lceil \frac{n}{2} \rceil+i-5]\}$. Thus, showing that $\textrm{Avg}(\mathcal{A}^{(6+\Delta)}) < \frac{n}{2}$ is sufficient for implying that $\textrm{Avg}(\mathcal{A}^{(k)}) < \frac{n}{2}$ for all such $k$. We compute $\textrm{Avg}(\mathcal{A}^{(6+\Delta)})$ to be $\frac{n^3+60n+16}{2n^2+4n+80}$ when $n$ is even and $\frac{n^3+3n^2+31n+29}{2n^2+8n+54}$ when $n$ is odd, in both cases less than $\frac{n}{2}$. 

\medskip

\noindent For $k \in \{2i+\Delta \ | \ i \in [3, \frac{n-\Delta-2}{2}]\}$, we have that $\mathcal{A}^{(k+1)} = \mathcal{A}^{(k)} \ \cup \ \{[\lceil \frac{n}{2} \rceil-\frac{k+2-\Delta}{2}]\}$ and $\mathcal{A}^{(k+2)} = \mathcal{A}^{(k)} \ \cup \ \{[\lceil \frac{n}{2} \rceil-\frac{k+2-\Delta}{2}], [\lceil \frac{n}{2} \rceil +\frac{k-5+\Delta}{2}]\}$. Therefore, $\mathcal{A}^{(k+1)} = \mathcal{A}^{(k)} \cup \ \mathcal{C}_1$ and $\mathcal{A}^{(k+2)} = \mathcal{A}^{(k)} \cup \ \mathcal{C}_2$, where $\mathcal{C}_1$ and $\mathcal{C}_2$ are families such that $\textrm{Avg}(\mathcal{C}_1) < \frac{n}{2}$ and $\textrm{Avg}(\mathcal{C}_2) < \frac{n}{2}$. Hence, $\textrm{Avg}(\mathcal{A}^{(k)}) < \frac{n}{2}$ implies that both $\textrm{Avg}(\mathcal{A}^{(k+1)}) < \frac{n}{2}$ and $\textrm{Avg}(\mathcal{A}^{(k+2)}) < \frac{n}{2}$. Noting that $\textrm{Avg}(\mathcal{A}^{(6+\Delta)}) < \frac{n}{2}$, it follows by induction that $\textrm{Avg}(\mathcal{A}^{(k)}) < \frac{n}{2}$ for all $k \in [7+\Delta, n]$. Finally, $\mathcal{A}^{(n+1)}=\mathcal{A}^{(n)} \cup \{\emptyset\}$ and $\textrm{Avg}(\mathcal{A}^{(n)}) < \frac{n}{2}$ together imply that $\textrm{Avg}(\mathcal{A}^{(n+1)}) < \frac{n}{2}$, completing the proof of Theorem 3.2.

\section*{4. Considering $3 \leq |\mathcal{B}| \leq 4$ for $h=4$}

It remains to consider the technique of averaging for when $3 \leq |\mathcal{B}| \leq 4$ (where $h$ is again equal to $4$). Extending Theorem 2.1 by these final two values of $|\mathcal{B}|$ would imply Conjecture 1.1 for $h=4$ (and more generally would conclude proof of $\textrm{Avg}(\mathcal{A}) \geq n/2$ for all separating union-closed families with $h=4$). We prove some respective propositions for $|\mathcal{B}| = 3$ and $|\mathcal{B}| = 4$.

\vspace{-0.125cm}

\subsection*{The case $|\mathcal{B}| = 3$}

We establish Propositions F, G, H, and I under the assumption that $|\mathcal{B}|=3$. Let $\mathcal{B} = \{B_1,B_2,B_3\}$, and for $i \in \{1,2,3\}$, let $k_i$ be the number of elements from $[n]$ that are contained in exactly $i$ member sets of $\mathcal{B}$.

\medskip

\noindent \textbf{Proposition F.} $|B| \in \{n-1, n\}$.

\medskip

\noindent \textit{Proof.} Otherwise, $|B| \leq n-2$, and $B_1 \subsetneq B_1 \cup B_2 \subsetneq B \subsetneq [n]$ contradicts Lemma 1.3.

\medskip

\noindent \textbf{Proposition G.} If $|B|=n$ and $A \in \mathcal{A}_{<n/2} \setminus \mathcal{B}$, then $\bigcup_{i=1}^3 \texttt{irr}_{\mathcal{B}}(B_i) \not \subseteq A$.

\medskip

\noindent \textit{Proof.} We let $|B|=n$, and assume that there exists $A \in \mathcal{A}_{<n/2} \setminus \mathcal{B}$ such that $\bigcup_{i=1}^3 \texttt{irr}_{\mathcal{B}}(B_i) \subseteq A$. By double counting, we have that $\sum_{i=1}^3 ik_i = \sum_{i=1}^3 |B_i|$. We further note that $k_1 +2(k_2+k_3) \leq \sum_{i=1}^3 ik_i$, $k_2+k_3=n-k_1$, and $|B_i| \leq \frac{n-1}{2}$ for every $i \in \{1,2,3\}$. Consequently, $k_1+2(n-k_1) \leq 3(\frac{n-1}{2})$, which implies that $k_1 \geq \frac{n+3}{2}$. Then $k_1 = \sum_{i=1}^3 |\texttt{irr}_{\mathcal{B}}(B_i)| = \ |\bigcup_{i=1}^3 \texttt{irr}_{\mathcal{B}}(B_i)|$ implies that $|A| \geq \frac{n+3}{2}$, which contradicts $A \in  \mathcal{A}_{<n/2}$. This concludes the proof of Proposition G.

\medskip
\smallskip

\noindent \textbf{Proposition H.} If $|B|=n$, then for all $A \in \mathcal{A}_{<n/2}$ and $i \in \{1,2,3\}$, $|\texttt{irr}_{\mathcal{B}}(B_i)| > 1$ implies that $|A \cap \texttt{irr}_{\mathcal{B}}(B_i)| \in \{0, |\texttt{irr}_{\mathcal{B}}(B_i)|-1, |\texttt{irr}_{\mathcal{B}}(B_i)|\}$.

\medskip

\noindent \textit{Proof.} Assume otherwise, i.e. that $|B|=n$ and there exists $A \in \mathcal{A}_{<n/2}$ and $i \in \{1,2,3\}$ such that both $|\texttt{irr}_{\mathcal{B}}(B_i)| > 2$ and $1 \leq |A \cap \texttt{irr}_{\mathcal{B}}(B_i)| \leq |\texttt{irr}_{\mathcal{B}}(B_i)|-2$. Without loss of generality, let $i=1$. Then $B_2 \subsetneq B_2 \cup B_3 \subsetneq  B_2 \cup B_3 \cup A \subsetneq [n]$ with $|B_2 \cup B_3 \cup A|<n-1$. This contradicts Lemma 1.3, completing the proof of Proposition H.

\medskip
\smallskip

\noindent \textbf{Proposition I.} If $|B|=n-1$ and $A \in \mathcal{A}_{<n/2} \setminus \mathcal{B}$, then either $A = \bigcup_{i=1}^3 \texttt{irr}_{\mathcal{B}}(B_i)$ or $A$ satisfies exactly one of the following three conditions:

\smallskip

\noindent (\romannumeral 1) $\hspace{0.0875cm} A \cap \texttt{irr}_{\mathcal{B}}(B_1) = \emptyset \ \land \ (B_2 \cup B_3) \setminus (B_2 \cap B_3) \subseteq A$;

\noindent (\romannumeral 2) $A \cap \texttt{irr}_{\mathcal{B}}(B_2) = \emptyset \ \land \ (B_1 \cup B_3) \setminus (B_1 \cap B_3) \subseteq A$;

\noindent (\romannumeral 3) $\hspace{-0.075cm} A \cap \texttt{irr}_{\mathcal{B}}(B_3) = \emptyset \ \land \ (B_1 \cup B_2) \setminus (B_1 \cap B_2) \subseteq A$.

\medskip

\noindent \textit{Proof.} We let $|B|=n-1$ and consider any $A \in \mathcal{A}_{n/2} \setminus \mathcal{B}$. For every $i \in \{1,2,3\}$, we have that $\texttt{irr}_{\mathcal{B}}(B_i) \cap A \in \{\emptyset, \texttt{irr}_{\mathcal{B}}(B_i)\}$. (Otherwise, without loss of generality $\emptyset \subsetneq \texttt{irr}_{\mathcal{B}}(B_1) \cap A \subsetneq \texttt{irr}_{\mathcal{B}}(B_1)$, which implies that $B_2 \subsetneq B_2 \cup B_3 \subsetneq B_2 \cup B_3 \cup A \subsetneq B_2 \cup B_3 \cup A \cup B_1 \subsetneq [n]$, contradicting $h=4$.)

\smallskip

\begin{itemize}

\item First, we assume that $\texttt{irr}_{\mathcal{B}}(B_i) \subseteq A$ for all $i \in \{1,2,3\}$. Applying the double counting argument from the proof of Proposition G to this case, we have that $k_1+2((n-1)-k_1) \leq 3(\frac{n-1}{2})$, which implies that $k_1 \geq \frac{n-1}{2}$. Noting that $k_1 = \sum_{i=1}^3 |\texttt{irr}_{\mathcal{B}}(B_i)| = \ |\bigcup_{i=1}^3 \texttt{irr}_{\mathcal{B}}(B_i)|$, we observe that if $A \setminus \bigcup_{i=1}^3 \texttt{irr}_{\mathcal{B}}(B_i) \neq \emptyset$, then $|A|>\frac{n}{2}$, which contradicts the fact that $A \in \mathcal{A}_{<n/2}$. Therefore, $A=\bigcup_{i=1}^3 \texttt{irr}_{\mathcal{B}}(B_i)$.

\smallskip

\item Next, we assume that that there is no $i \in \{1,2,3\}$ such that $\texttt{irr}_{\mathcal{B}}(B_i) \subseteq A$. In this case, $A \subsetneq A \cup B_1 \subsetneq A \cup B_1 \cup B_2 \subsetneq A \cup B_1 \cup B_2 \cup B_3 \subsetneq [n]$, contradicting $h=4$.

\smallskip

\item Now, we assume that that there exists exactly one $i \in \{1,2,3\}$ such that $\texttt{irr}_{\mathcal{B}}(B_i) \subseteq A$. Without loss of generality, let $i=1$. Then $(B_1 \subsetneq A \cup B_1 \subsetneq A \cup B_1 \cup B_2 \subsetneq B \subsetneq [n]) \ \lor \ (A \subsetneq B_1 \subsetneq B_1 \cup B_2 \subsetneq B \subsetneq [n])$, again a contradiction with $h=4$.

\smallskip

\item Finally, we assume that there exist exactly two distinct elements $i$ and $j$ in $\{1,2,3\}$ such that $\texttt{irr}_{\mathcal{B}}(B_i) \subseteq A$ and $\texttt{irr}_{\mathcal{B}}(B_j) \subseteq A$. Denote by $k$ the unique element from the set $\{1,2,3\}\setminus\{i,j\}$. If $(B_i \cap B_k) \setminus B_j \not \subseteq A$, then $B_j \subsetneq B_j \cup A \subsetneq B_i \cup B_j \cup A \subsetneq B \subsetneq [n]$. Similarly, if $(B_j \cap B_k) \setminus B_i \not \subseteq A$, then $B_i \subsetneq B_i \cup A \subsetneq B_i \cup B_j \cup A \subsetneq B \subsetneq [n]$. In either case, there exists a chain of size $5$ in $\mathcal{A}$, contradicting $h=4$. It follows that $(B_i \cap B_k) \setminus B_j \subseteq A$ and $(B_j \cap B_k) \setminus B_i \subseteq A$ (as well as $\texttt{irr}_{\mathcal{B}}(B_i) \subseteq A$ and $\texttt{irr}_{\mathcal{B}}(B_j) \subseteq A$), as illustrated in Figure 4.1. 

\end{itemize}

\medskip

\noindent Hence, if $A \neq \bigcup_{i=1}^3 \texttt{irr}_{\mathcal{B}}(B_i)$, then $A$ satisfies exactly one of (\romannumeral 1), (\romannumeral 2), or (\romannumeral 3). This completes the proof of Proposition I.

\medskip

\def\firstcircle{(90:1.75cm) circle (2.5cm)}
\def\secondcircle{(210:1.75cm) circle (2.5cm)}
\def\thirdcircle{(330:1.75cm) circle (2.5cm)}
  
\begin{figure}[H]
\captionsetup{width=0.89\linewidth}
\caption*{\textbf{Figure 4.1:} If $|\mathcal{B}|=3$ and $|B|=n-1$, then any $A \in \mathcal{A}_{<n/2}\setminus \mathcal{B}$ such that $A \neq \bigcup_{i=1}^3 \texttt{irr}_{\mathcal{B}}(B_i)$ must have one of three forms.}

\vspace{-0.375cm}

\begin{center}
      
\begin{tikzpicture}[scale=0.6375]

        \node at (-0.0275,4.875){\large{(\romannumeral 1)}};
        \begin{scope}
        \fill[lightgray] \secondcircle;
        \fill[lightgray] \thirdcircle;
        \end{scope}
        \begin{scope}
        \clip \secondcircle;
        \fill[white] \thirdcircle;
        \end{scope}
        \draw \secondcircle;
        \node at (0,2.3){$\texttt{irr}_{\mathcal{B}}(B_1)$};
        \draw \thirdcircle;
        \node at (-2.45,-1.25) {$\texttt{irr}_{\mathcal{B}}(B_2)$};
        \draw \firstcircle;
        \node at (2.45,-1.25){$\texttt{irr}_{\mathcal{B}}(B_3)$};
        \node at (0,0.285){\LARGE{$\therefore$}};
        \node[rotate=180] at (0,-1.875){\LARGE{$\therefore$}};
        \node at (0,3.5){\underline{$B_1$}};
        \node at (-1.75,-2.5){\underline{$B_2$}};
        \node at (1.75,-2.5){\underline{$B_3$}};
    
\end{tikzpicture}
    
\vspace{-0.25cm}
    
\begin{tikzpicture}[scale=0.6375]

        \node at (0,4.875){\large{(\romannumeral 2)}};
        \begin{scope}
        \fill[lightgray] \thirdcircle;
        \fill[lightgray] \firstcircle;
        \end{scope}
        \begin{scope}
        \clip \thirdcircle;
        \fill[white] \firstcircle;
        \end{scope}
        \draw \thirdcircle;
        \node at (0,2.3){$\texttt{irr}_{\mathcal{B}}(B_1)$};
        \draw \firstcircle;
        \node at (-2.45,-1.25) {$\texttt{irr}_{\mathcal{B}}(B_2)$};
        \draw \secondcircle;
        \node at (2.45,-1.25){$\texttt{irr}_{\mathcal{B}}(B_3)$};
        \node at (0,0.285){\LARGE{$\therefore$}};
        \node[rotate=65] at (1.05,0.865){\LARGE{$\therefore$}};
        \node at (0,3.5){\underline{$B_1$}};
        \node at (-1.75,-2.5){\underline{$B_2$}};
        \node at (1.75,-2.5){\underline{$B_3$}};
        
\end{tikzpicture}
\hspace{0.85cm}
\begin{tikzpicture}[scale=0.6375]
        
        \node at (-0.025,4.875){\large{(\romannumeral 3)}};
        \begin{scope}
        \fill[lightgray] \firstcircle;
        \fill[lightgray] \secondcircle;
        \end{scope}
        \begin{scope}
        \clip \firstcircle;
        \fill[white] \secondcircle;
        \end{scope}
        \draw \firstcircle;
        \node at (0,2.3){$\texttt{irr}_{\mathcal{B}}(B_1)$};
        \draw \secondcircle;
        \node at (-2.45,-1.25) {$\texttt{irr}_{\mathcal{B}}(B_2)$};
        \draw \thirdcircle;
        \node at (2.45,-1.25){$\texttt{irr}_{\mathcal{B}}(B_3)$};
        \node at (0,0.285){\LARGE{$\therefore$}};
        \node[rotate=57] at (-1.535,0.85){\LARGE{$\therefore$}};
        \node at (0,3.5){\underline{$B_1$}};
        \node at (-1.75,-2.5){\underline{$B_2$}};
        \node at (1.75,-2.5){\underline{$B_3$}};

\end{tikzpicture}
        
\bigskip
\medskip  
        
\hspace{-3cm} $^* \ $ Gray indicates that the region is a subset of $A$.
    
\hspace{-0.2cm} $^{**}$ Three dots indicate that part of the region may be a subset of $A$.
    
\end{center}

\end{figure}
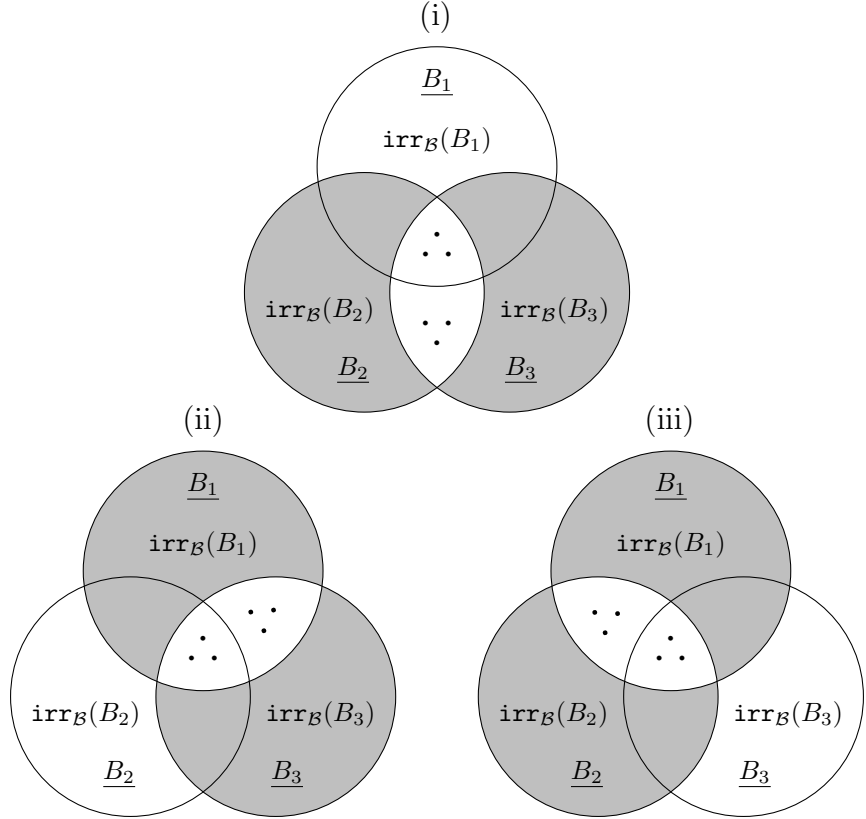

\subsection*{The case $|\mathcal{B}| = 4$}

We conclude this study with propositions valid for $|\mathcal{B}|=4$, culminating in Theorem 4.1. 

\medskip

\noindent Let $\mathcal{B}=\{B_1,B_2,B_3,B_4\}$.

\medskip

\noindent \textbf{Proposition J.} $|B|=n$.

\medskip

\noindent \textit{Proof.} If not, then $B_1 \subsetneq B_1 \cup B_2 \subsetneq B_1 \cup B_2 \cup B_3 \subsetneq B \subsetneq [n]$, which contradicts $h=4$.

\medskip

\noindent \textbf{Proposition K.} For each $i \in \{1,2,3,4\}$, $|\texttt{irr}_{\mathcal{B}}(B_i)|=1$.

\medskip

\noindent \textit{Proof.} Otherwise, let $|\texttt{irr}_{\mathcal{B}}(B_1)| > 1$ without loss of generality. It follows that $B_2 \subsetneq B_2 \cup B_3 \subsetneq B_2 \cup B_3 \cup B_4 \subsetneq [n]$ with $|B_2 \cup B_3 \cup B_4| < n-1$, contradicting Lemma 1.3.

\medskip

\noindent \textbf{Proposition L.} (\romannumeral 1) If $n$ is even, then $\sum_{i=1}^4 |B_i| = 2n-4$ and $|B_i|=\frac{n-2}{2}$ for all $i \in \{1,2,3,4\}$. (\romannumeral 2) If $n$ is odd, then $2n-4 \leq \sum_{i=1}^4 |B_i| \leq 2n-2$ and $|B_i| \geq \frac{n-5}{2}$ for all $i \in \{1,2,3,4\}$. Further, if there exists $i \in \{1,2,3,4\}$ such that $|B_i|=\frac{n-5}{2}$, then $|B_j|=\frac{n-1}{2}$ for all $j \in \{1,2,3,4\} \setminus \{i\}$.

\medskip

\noindent \textit{Proof.} For $i \in \{1,2,3,4\}$, let $k_i$ be the number of elements from $[n]$ that are contained in exactly $i$ member sets of $\mathcal{B}$. We observe that $k_1=4$ by Proposition K, and that $n-k_1 = k_2+k_3+k_4$. Additionally, double counting gives $\sum_{i=1}^4 ik_i = \sum_{i=1}^4 |B_i|$. Therefore, we have that $k_1+2(k_2+k_3+k_4) = 2n-4 \leq \sum_{i=1}^4 |B_i|$.

\medskip

\noindent For part (\romannumeral 1), we assume that $n$ is even. Then $\sum_{i=1}^4 |B_i| \leq 4(\frac{n-2}{2})=2n-4$, which together with the above inequality $2n-4 \leq \sum_{i=1}^4 |B_i|$ implies that $\sum_{i=1}^4 |B_i| = 2n-4$. Because $n$ is even, $|B_i| \leq \frac{n-2}{2}$ for every $i \in \{1,2,3,4\}$. Together with $\sum_{i=1}^4 |B_i| = 2n-4$, this then implies that $|B_i|=\frac{n-2}{2}$ for every $i \in \{1,2,3,4\}$.

\medskip

\noindent For part (\romannumeral 2), we assume that $n$ is odd. Then $\sum_{i=1}^4 |B_i| \leq 4(\frac{n-1}{2})=2n-2$, so in this case we have that $2n-4 \leq \sum_{i=1}^4 |B_i| \leq 2n-2$. It follows that for every $i \in \{1,2,3,4\}$, $|B_i| \geq \frac{n-5}{2}$. If there exists $i \in \{1,2,3,4\}$ such that $|B_i|=\frac{n-5}{2}$, then $\sum_{j \in \{1,2,3,4\} \setminus \{i\}}|B_j| \geq 2n-4 - \frac{n-5}{2} = \frac{3n-3}{2}$. Since $B_j \in \mathcal{A}_{< n/2}$ for all $j \in \{1,2,3,4\} \setminus \{i\}$, we then have that $|B_j| = \frac{n-1}{2}$ for all such $j$, completing the proof of Proposition L.

\medskip
\medskip

\noindent \textbf{Theorem 4.1.} \textit{For any separating union-closed family} $\mathcal{A}$ \textit{with} $h=|\mathcal{B}|=4$:

\[\textrm{Avg}(\mathcal{A}) \ > \ \Bigr \lfloor \frac{n}{2} \Bigr \rfloor - 1\textrm{.}\]

\medskip
\medskip

\noindent \textit{Proof.} The following blocks comprise the proof of this theorem:

\smallskip

\begin{itemize}

\item[1.] For each $i \in \{1,2,3,4\}$, we denote by $b_i$ the unique element from $\texttt{irr}_{\mathcal{B}}(B_i)$ (see Proposition K). We consider any $A \in \mathcal{A}_{<n/2} \setminus \mathcal{B}$, and note that $b_i \in A$ for some $i \in \{1,2,3,4\}$. (Otherwise, $A \subsetneq A \cup B_1 \subsetneq A \cup B_1 \cup B_2 \subsetneq A \cup \bigcup_{i=1}^3 B_i \subsetneq [n]$, contradicting $h=4$.) Without loss of generality, we assume that $b_1 \in A$ and let $\mathcal{B}'=\{A,B_2,B_3,B_4\}$. Noting that $b(\mathcal{B}) = b(\mathcal{B}')$, it must be that $\mathcal{B}'$ is an irredundant subfamily of $\mathcal{A}$. (If not, then $|\mathcal{B}|<4$, a contradiction.) Then, because $\mathcal{B}$ is defined as any smallest irredundant subfamily of $\mathcal{A}_{<n/2}$ such that $b(\mathcal{B})=b(\mathcal{A}_{<n/2})$, and $\mathcal{B'}$ satisfies these conditions, the propositions of this section also apply to $\mathcal{B}'$.

\medskip

\item[2.] We assume that $n$ is even. By Proposition L, $|B'| = \frac{n-2}{2}$ for every $B' \in \mathcal{B}'$, implying that $|A|=\frac{n-2}{2}$. Therefore, any member set from $\mathcal{A}_{<n/2}$ must have size equal to $\frac{n-2}{2}$, implying that every member set in $\mathcal{A}$ has size greater than or equal to $\frac{n-2}{2}$. Noting that $[n] \in \mathcal{A}$, we then have the result that $\textrm{Avg}(\mathcal{A}) > \frac{n-2}{2}=\lfloor \frac{n}{2} \rfloor - 1$.

\medskip

\item[3.] We now assume that $n$ is odd. 

\smallskip

\begin{itemize}

\item[A.] There must be some element $j \in \{2,3,4\}$ such that $b_j \in A$. (Otherwise, we again have a contradiction with $h=4$ in that $A \subsetneq A \cup B_2 \subsetneq A \cup B_2 \cup B_3 \subsetneq A \cup \bigcup_{i=2}^4 B_i \subsetneq [n]$.) Without loss of generality, we assume that $b_2 \in A$ and let $\mathcal{B}''=\{B_1,A,B_3,B_4\}$. $\mathcal{B}''$ is yet another family satisfying the conditions of $\mathcal{B}$, and is thus subject to the propositions of this section. We assume that there exists $i \in \{1,2,3,4\}$ such that $|B_i|=\frac{n-5}{2}$. Recall that, by part (\romannumeral 2) of Proposition L, if there exists $i \in \{1,2,3,4\}$ such that $|B_i|=\frac{n-5}{2}$, then $|B_j|=\frac{n-1}{2}$ for all $j \in \{1,2,3,4\} \setminus \{i\}$. Therefore, we have that $|A|=\frac{n-1}{2}$ by applying Proposition L to $\mathcal{B}''$ if $|B_1|=\frac{n-5}{2}$, to $\mathcal{B}'$ if $|B_2|=\frac{n-5}{2}$, and to $\mathcal{B}'$ (or $\mathcal{B}''$) if $|B_3|=\frac{n-5}{2}$ or $|B_4|=\frac{n-5}{2}$. We now assume that no member set of $\mathcal{B}$ has size $\frac{n-5}{2}$, i.e. that $|B_i|\geq \frac{n-3}{2}$ for every $i \in \{1,2,3,4\}$. If $|A| \leq \frac{n-5}{2}$, then $|A|=\frac{n-5}{2}$ by Proposition L, and we further apply Proposition L to $\mathcal{B}''$ in order to obtain that $|B_1|=\frac{n-1}{2}$, to $\mathcal{B}'$ in order to obtain that $|B_2|=\frac{n-1}{2}$, and to $\mathcal{B}'$ (or $\mathcal{B}''$) in order to obtain that $|B_3|=\frac{n-1}{2}$ and $|B_4|=\frac{n-1}{2}$.

\medskip

\item[B.] We execute block 1 of this proof followed by block 3A, applying the local reassignment $A \leftarrow A'$ and $\mathcal{B} \leftarrow \mathcal{B}'$, in order to obtain that $|A'|=\frac{n-1}{2}$. $A$ is therefore the only member set of $\mathcal{A}$ with size less than or equal to $\frac{n-5}{2}$. Again noting that $[n] \in \mathcal{A}$, this implies that $\textrm{Avg}(\mathcal{A}) > \frac{n-3}{2} = \lfloor \frac{n}{2} \rfloor - 1$, completing the proof of Theorem 4.1.

\end{itemize} 

\end{itemize}

\end{document}